\documentclass[11pt]{amsart}
\usepackage[colorlinks=true, pdfstartview=FitV, linkcolor=blue, citecolor=blue, urlcolor=blue]{hyperref}

\newtheorem{theorem}{Theorem}[section]
\newtheorem{lemma}[theorem]{Lemma}
\newtheorem{corollary}[theorem]{Corollary}
\newtheorem{proposition}[theorem]{Proposition}
\newtheorem{remark}[theorem]{Remark}

\def\grad{\operatorname{grad}}
\def\int{\operatorname{int}}
\def\bd{\operatorname{bd}}

\def\codim{\operatorname{codim}}

\def\Ad{\operatorname{Ad}}
\def\ad{\operatorname{ad}}

\def\lie#1{\mathfrak{ #1}}

\def\liek{\lie k}
\def\liep{\lie p}

\def\i{\mathrm{i}}
\def\m{\mathcal{M}}
\def\mp{\mathcal{M}_{\lie p}}

\newcommand{\ip}{{\im\mathfrak p}}
\newcommand{\lieg}{{\mathfrak g}}
\newcommand{\liea}{{\mathfrak a}}
\begin{document}

%%%%%%%%%%%%%%%%%%%%%%%%%%%%%%%%%%%%%%%%%%%%%%%%%%%%%%%%%%%%%%%%%%%
\newcommand{\acts}{\,\mbox{\raisebox{0.27ex}{\tiny{$\bullet$}}}\,}
\newcommand{\im}{\mathrm{i}}
\newcommand{\sacts}{\,\mbox{\raisebox{-0.5ex}{\Large{$\cdot$}}}}
\newcommand{\ssacts}{\,\mbox{\raisebox{-0.4ex}{\Large{$\cdot$}}}}
\newcommand{\sodass}{\ \vert \ }
\newcommand{\inv}{^{-1}}
\newcommand{\seq}[1]{({#1})_{n \in \nn}}
\newcommand{\sperp}{{\mbox{\tiny {$\perp$}}}}
%%%%%%%%%%%%%%%%%%%%%%%%%%%%%%%%%%%%%%%%%%%%%%%%%%%%%%%%%%%%%%%%%%
\newcommand{\mgc}{{\mg}^{\mbox{\tiny $\cc$}}}
\newcommand{\munc}{\mun^{\mbox{\tiny $\cc$}}}
\newcommand{\kc}{K^\mathbb C}
\newcommand{\gx}{{G \acts x}}
\newcommand{\gv}{{G \acts v}}
\newcommand{\kbeta}{{K \acts \beta}}
\newcommand{\obeta}{\mathcal O_\beta}
%%%%%%%%%%%%%%%%%%%%%%%%%%%%%%%%%%%%%%%%%%%%%%%%%%%%%%%%%%%%%%%%%

%%
\newcommand{\ia}{\mathfrak a_+^\ast}
\newcommand{\iao}{\mathfrak a_0^\ast}
\newcommand{\tma}{\tilde {\mathfrak a}}
%%%%%%%%%%%%%%%%%%%%%%%%%%%%%%%%%%%%%%%%%%%%%%%%%%%%%%%%%%%%%%%%%%%
\newcommand{\mip}{{\mu_{\mathfrak p}}}
\newcommand{\muip}{{\mu_{\im \mathfrak p}}}
\newcommand{\mup}{\mu_{\mathfrak p}}
\newcommand{\mia}{\mu_{\mathfrak a}}
\newcommand{\mua}{\mu_{\mathfrak a}}
\newcommand{\mipa}{{\mu_{\ip,\alpha}}}
\newcommand{\mipb}{{\mu_{\lie p,\beta}}}
\newcommand{\dmip}{\mathrm d \mip(x)}
\newcommand{\mipn}{\mu_\lie p^{(n)}}
%%%%%%%%%%%%%%%%%%%%%%%%%%%%%%%%%%%%%%%%%%%%%%%%%%%%%%%%%%%%%%%%%%%
\newcommand{\Mip}{\mathcal M_{\liep}}
\newcommand{\Mipb}{\mathcal M_{\lie p,\beta}}
\newcommand{\Mipa}{\mathcal M_{\lie p,\alpha}}
\newcommand{\tMip}{\tilde {\mathcal M}_\lie p}
\newcommand{\SGMip}[1]{\mathcal S_{G}(\Mip({#1}))}
\newcommand{\SGMipb}[1]{\mathcal S_{G}(\Mipb({#1}))}
\newcommand{\tSGMip}{\mathcal S_{G}(\tilde{\mathcal M}_\ip)}
\newcommand{\SGM}{\mathcal S_{G}(\Mip)}
\newcommand{\SAM}{\mathcal S_A(\mathcal M_{i \ma})}
%%%%%%%%%%%%%%%%%%%%%%%%%%%%%%%%%%%%%%%%%%%%%%%%%%%%%%%%%%%%%%%%%%%
\newcommand{\zcone}{\mathcal N_{G}}
\newcommand{\zconess}{\mathcal N_{G_s}}
\newcommand{\Vohne}{V^\circ}
%%%%%%%%%%%%%%%%%%%%%%%%%%%%%%%%%%%%%%%%%%%%%%%%%%%%%%%%%%%%%%%%%%%
\newcommand{\xa}{X_{\alpha}}
\newcommand{\xb}{X_{\beta}}
\newcommand{\xc}{X_{(T_\gamma)}}
\newcommand{\Xa}{X^{T_\alpha}}
%%%%%%%%%%%%%%%%%%%%%%%%%%%%%%%%%%%%%%%%%%%%%%%%%%%%%%%%%%%%%%%%%%%%
\newcommand{\pn}{{\mathbb P^n(\mathbb C)}}
\newcommand{\hyp}[1]{\mathcal O_{#1}} 
%%%%%%%%%%%%%%%%%%%%%%%%%%%%%%%%%%%%%%%%%%%%%%%%%%%%%%%%%%%%%%%%%%
\newcommand{\zs}{Z_\beta^{\mbox{\raisebox{0.3ex}{\scriptsize $ss$}}}}
\newcommand{\zss}{Z_{\beta}^{\mbox{\raisebox{0.3ex}{\scriptsize $ss$}}} {\mbox{\raisebox{0.2ex} {\scriptsize$(+)$}}}}
\newcommand{\zplus}{Z_\beta{\mbox{\raisebox{0.2ex} {\scriptsize$(+)$}}}}
%%%%%%%%%%%%%%%%%%%%%%%%%%%%%%%%%%%%%%%%%%%%%%%%%%%%%%%%%%%%%%%%%%%
\def\twist#1#2#3{#1\times^{#2}#3}

\newcommand{\ma}{\mathfrak a}
\newcommand{\mb}{\mathfrak b}
\newcommand{\mg}{\mathfrak g}
\newcommand{\mh}{\mathfrak h}
\newcommand{\mk}{\mathfrak k}
\newcommand{\ml}{\mathfrak l}
\newcommand{\mm}{\mathfrak m}
\newcommand{\mn}{\mathfrak n}
\newcommand{\mq}{\mathfrak q}
\newcommand{\mr}{\mathfrak r}
\newcommand{\ms}{\mathfrak s}
\newcommand{\mt}{\mathfrak t}
\newcommand{\mun}{\mathfrak u}
\newcommand{\lieu}{\mathfrak u}
\newcommand{\mz}{\mathfrak z}
%%%%%%%%%%%%%%%%%%%%%%%%%%%%%%%%%%%%%%%%%%%%%%%%%%%%%%%%%%%%%%%%%%%%
\newcommand{\cc}{\mathbb C}
\newcommand{\CC}{\mathbb C}

\newcommand{\rr}{\mathbb R}
\newcommand{\RR}{\mathbb R}
\newcommand{\zz}{\mathbb Z}
\newcommand{\nn}{\mathbb N}
\newcommand{\pp}{\mathbb P}
\newcommand{\PP}{\mathbb P}
%%%%%%%%%%%%%%%%%%%%%%%%%%%%%%%%%%%%%%%%%%%%%%%%%%%%%%%%%%%%%%%%%%%
\newcommand{\qd}{\hfill{\( \square \)}}
%%%%%%%%%%%%%%%%%%%%%%%%%%%%%%%%%%%%%%%%%%%%%%%%%%%%%%%%%%%%%%%%%%%%
\newcommand{\bn}{\bigskip \noindent}
\newcommand{\n}{\noindent}
\newfont{\mib}{cmmib10}

%%%%%%%%%%%%%%%%%%%%%%%%%%%%%%%%%%%%%%%%%%%%%%%%%%%%%%%%%%%%%%%%%

\title[Convexity of gradient maps]{Convexity properties of gradient maps}
\author{Peter Heinzner and Patrick Sch\"utzdeller}
\address{Fakult\"at und Institut f\"ur Mathematik \\ 
               Ruhr-Universit\"at Bochum \\
               D-44780 Bochum}
\email{heinzner@cplx.rub.de}
\address{Institut f\"ur Mathematik, Universit\"at Paderborn, D-33095 Paderborn}
\email{schuetzd@math.upb.de}
\thanks{Both authors are partially supported by the Sonderforschungsbereich SFB TR 12 of the Deutsche Forschungsgemeinschaft} 
\date{}
\maketitle

\begin{abstract}
We consider the action of a real reductive group $G$ on a K\"ahler manifold $Z$ which is the restriction of a holomorphic action of the complexified group $G^\cc$. We assume that the induced action of a compatible maximal compact subgroup $U$ of $G^\cc$ on $Z$ is Hamiltonian. We have an associated gradient map $\mip \colon Z \to \liep$ where $ \lieg = \liek \oplus \liep$ is the Cartan decomposition of $\lie g$.
For a $G$-stable subset $Y$ of $Z$ we consider convexity properties of the intersection of $\mip(Y)$ with a closed Weyl chamber in a maximal abelian subspace $\lie a$ of $\liep$. Our main result is a Convexity Theorem for real semi-algebraic subsets $Y$ of $Z = \pp(V)$ where $V$ is a unitary representation of $U$.
\end{abstract}

\section{Introduction}
Let $U$ be a compact Lie group and $U^\CC$ its complexification. Then the map
map $U\times \im \lie u\to U^\CC$, $(u,\xi)\mapsto u\exp\xi$ is a diffeomorphism.
A closed subgroup $G$ of $U^\CC$ with Lie algebra $\lieg$ is said to be compatible
if the restriction $K\times \liep\to G$ is a diffeomorphism where $K=G\cap U$
and $\liep=\im\lieu\cap \lie g$. In the rest of this paper we fix a compatible $G$
a compact complex manifold $Z$ and a holomorphic action $U^\CC\times Z\to Z$.
We also assume that there is a $U$-invariant K\"ahler form $\omega$ and a $U$-equivariant
momentum map $\mu\colon Z \to \lieu^*$. We fix a $U$-invariant inner product $\langle , \rangle$ on 
$\lie u\cong \im\lieu$ and view $\mu$ as a map from $Z$ into $\im \lieu$.
Since $\liep\subset\im\lieu$ the composition of $\mu$ with the orthogonal projection
of $\im\lieu$ onto $\liep$ defines a $K$-equivariant map $\mup\colon Z\to \liep$
which we call the $G$-gradient map. Then we have $\grad\mup^\xi=\xi_Z$
for $\xi\in\liep$
where $\grad$ is computed with respect to the Riemannian structure given 
by $\omega$, $\mup^\xi:=\langle \mup,\xi\rangle$ and $\xi_Z$ is the vector field
induced by the action. For a maximal dimensional Lie subalgebra $\lie a$ of $\lie g$ 
which is contained in $\lie p$ and a $G$-stable subset $Y$ of $Z$ we have the
set $A(Y):=\mup(Y)\cap\liea$. In section \ref{semistable} we prove the following

\medskip
\noindent
{\bf Theorem.} \ \it
If $Y$ is closed, then $A(Y)$ is a finite union of convex polytopes. Each of the
polytopes is the convex hull of $\mup$-images of $A$ fixed points in $Y$
where $A=\exp\liea$. \rm  

\medskip
Let $\liea_+$ be a positive Weyl chamber of $\liea$ and set 
$A(Y)_+:=A(Y)\cap \liea_+$. The main 
result of this paper is the following 

\medskip
\noindent
{\bf Convexity Theorem.} \ \it
Let $Z=\PP(V)$ be the complex projective space of a unitary representation of $U$ 
and $\omega$ the induced K\"ahler structure on $\PP(V)$ 
with standard momentum map. Let $Y$ be a closed real semi-algebraic
subset of $Z$ whose real algebraic Zariski closure is irreducible. Then
$A(Y)_+$ is a convex polytope.
\rm

\medskip
\noindent
{\bf Corollary.} \ \it
Let $Z = U^\cc/Q$ be a complex flag manifold endowed with a $U$-invariant 
K\"ahler structure and $G$ a real form of $U^\CC$. 
Then $A(\overline{G \cdot x})_+$ is a convex polytope for any $x\in Z$.
\rm

\medskip
All of the above are well known if $G$ is a complex subgroup of $U^\CC$ 
and $Y$ is an irreducible complex analytic subset of $Z$. In this entirely 
holomorphic setup the Convexity Theorem holds for any compact K\"ahler
manifold $Z$ (see e.g. \cite{GuilleminSjamaar} for more on the history 
of the subject). On the other hand very little is known for a non 
complex group $G$ and general $Z$. In this generality convexity 
of $A(Y)_+$ is known only in very special cases. See e.\,g. \cite{Kostant} and \cite{OS}

\section{Basic properties 
of the gradient map}
\label{section:gradientmaps}
As before we assume that $U^\CC$ acts holomorphically on a compact
K\"ahler manifold $Z$ and that the K\"ahler form $\omega$ is
$U$-invariant. It is also assumed that there is a $U$-equivariant
momentum map $\mu$ and we denote the associated $G$-gradient map 
by $\mup$. For the convenience of the reader we recall here several
known basic facts which will be needed later.  

For a subspace $\lie m$ of $\lieg $ and $z\in Z$ let $\lie m \cdot
z:=\{\xi_Z(z) \sodass \xi\in \lie m\}$. The following elementary fact is
shown in \cite{HeinznerSchwarz}.

\begin{lemma}\label{lemma:kerdmu} 
  We have $\ker d\mup(z)=(\liep \cdot z)^\perp$ for all $z\in Z$.
\end{lemma}

Let $G=K\exp\liep$ be a compatible closed subgroup of $U^\CC$.  For
$\beta\in\liep$ we set $\mp(\beta):=\mup\inv(\beta)\subset Z$ and 
$\mp:=\mp(0)$.  For $z\in \mp$ the isotropy group $G_z=K_z\exp\liep_z$
is a compatible subgroup of $U^\CC$ (\cite[5.5]{HeinznerSchwarz}).
Since the $G_z$-representation on $T_z(Z)$ is completely reducible
(\cite[14.9]{HeinznerSchwarz}), there is a $G_z$-stable decomposition
$T_z(Z)=\lieg\cdot z\oplus W$. We have the following general Slice
Theorem (\cite[14.10, 14.21]{HeinznerSchwarz}):

\begin{theorem}[Slice Theorem] \label{slicetheorem} Let $z\in \mp$.
  Then there exists a $G_z$-stable open neighborhood $S$ of $0\in W$,
  a $G$-stable open neighborhood $\Omega$ of $z\in Z$ and a
  $G$-equivariant diffeomorphism $\Psi\colon \twist G{G_z}S\to \Omega$
  where $\Psi([e,0])=z$ and $\twist G{G_z}S$ denotes the $G$-bundle associated 
  with the principal bundle $G\to G/G_z$.
\end{theorem}

Actually, we have a Slice Theorem at every $z\in Z$. Set
$\beta:=\mup(z)$ and let $G^\beta=\{g\in G: \Ad g\cdot \beta=\beta\}$
denote the centralizer of $\beta$. Then we have a slice for the action
of $G^\beta$, as follows.

The centralizer $G^\beta$ is a compatible subgroup of $U^\CC$ with
Cartan decomposition $G^\beta=K^\beta\exp(\liep^\beta)$ where
$K^\beta=K\cap G^\beta$ and $\liep^\beta=\{\xi\in\liep:
\ad(\xi)\beta=0\}$. The group $G^\beta$ is also compatible with the
Cartan decomposition of the centralizer $(U^\cc)^\beta=(U^\beta)^\CC$ and $\beta$ is
fixed by the action of $U^\beta$ on $\lieu^\beta$. This implies that
the $\lieu^\beta$-component of $\mu$ defines a $U^\beta$-equivariant
shifted momentum map $\widehat{\mu_{\lieu^\beta}}\colon Z\to
\lieu^\beta$,
$\widehat{\mu_{\lieu^\beta}}(z)=\mu_{\lieu_\beta}(z)-\beta$. The
associated $G^\beta$-gradient map is given by
$\widehat{\mu_{\liep^\beta}}\colon Z\to \liep^\beta$,
$\widehat{\mu_{\liep^\beta}}(z)=\mu_{\liep^\beta}(z)-\beta$. This
shows that the Slice Theorem applies to the action of $G^\beta$ at
every point $z\in
(\widehat{\mu_{\liep^\beta}})\inv(0)=\mathcal{M}_{\liep^\beta}(\beta)$.
In particular, if $G$ is commutative, then we have a Slice Theorem for
$G$ at every point of $Z$.

\section{Orbit-type stratification}
Let $Y$ be a closed $G$-stable subset of $Z$. Our initial goal is to 
show that $A(Y) = \mip(Y) \cap \liea$ is a finite union of convex polytopes.
The proof is completed in section~\ref{semistable}. It depends on the orbit 
type stratification of $Z$ with respect to compatible commutative subgroups of $G$
which we explain now. 

Let $\lie b$ be a Lie subalgebra of $\lieg$ which is contained in
$\liep$. Note that $B=\exp \lie b$ is compatible with the Cartan
decomposition of $U^\CC$ and that $\exp\colon \lie b\to B$ is an
isomorphism of commutative Lie groups. Let $Z^B$ denote the set of
$B$-fixed points in $Z$. Since $B=\exp \lie b$ we have $Z^B=Z^{\lie
  b}=\{x\in Z: \xi_Z(x)=0 \text{ for all } \xi\in\lie b\}$. Further let 
$\mu_\lie b \colon Z \to \lie b$ the composition of $\mu$ with the orthogonal projection of $\im \lie u$ onto $\lie b$.

\begin{lemma} \label{lemma:imageoffixedpoints} For $B=\exp \lie b$ we
  have
  \begin{enumerate}
  \item the set $Z^B$ of $B$-fixed points is a smooth complex 
  submanifold of $Z$ and 
  \item the $B$-gradient
  map $\mu_{\lie b}\vert_{Z^B}\colon Z^B\to \lie b$
    is locally constant.
  \end{enumerate}
\end{lemma}
\begin{proof}
  Since $B$ acts on $Z$ by holomorphic transformations the set $Z^B$
  is a complex subspace of $Z$. The isotropy representation defines a
  linear $B$-action on $T_x(Z)$. By the Slice Theorem a $B$-stable
  open neighborhood of $x$ is $B$-equivariantly diffeomorphic to an
  open neighborhood of $0$ in $T_x(Z)$. Since the set of fixed points
  of a linear action is a linear subspace the set $Z^B$ is smooth.
  This shows (1).  The second assertion follows from
  Lemma~\ref{lemma:kerdmu}.
\end{proof}

For a connected subgroup $B$ of $A = \exp(\lie a)$ let $Z^{(B)}:=\{z\in Z: A_z=B\}$.  
The group $B$ is compatible and we have $Z^{(B)}=Z^{(\lie b)}:=\{z\in Z:
\liea_z=\lie b\}$.  A connected component $S$ of $Z^{(B)}$ is called
an $A$-stratum of type $A/B$ or alternatively an $\liea$-stratum
of type $\liea/\lie b$.

\begin{lemma} Let $S$ be an $\liea$-stratum of type $\liea/\lie
  b$, $q\in\mua(S)$ and $\liea(S):=q+\lie b^\perp$. Then we have:
  \begin{enumerate}
  \item $S$ is open in $Z^B$.
  \item $\mua(S)$ is an open subset of $\liea(S)$.
  \item $\mua\colon S\to \liea(S)$ is a submersion.
  \end{enumerate}
\end{lemma}
\begin{proof}
  Any $x\in Z^{(B)}$ has an open $A$-stable neighborhood $\Omega$
  which is $A$-equivariantly diffeomorphic to an $A$-stable
  neighborhood of $[e,0]$ in $\twist A B W$ where $W$ is a
  $B$-representation space and $e$ is the neutral element in $A$ 
  (Slice Theorem).  The $A$-stratum in
  $\twist A B W$ of type $A/B$ is given by $\twist A B W^B$ and
  coincides with the set of $B$-fixed points in $\twist A B W$. This
  shows that $Z^{(B)}$ is open in $Z^B$. In particular $S$ is open in
  $Z^B$ and we have (1).

  With respect to the orthogonal decomposition $\liea=\lie b\oplus
  \lie b^\perp$ we have $q=q_\lie b+ q_{\lie b^\perp}\in \mua(S)$ and
  $\mua=\mu_{\lie b}\oplus \mu_{\lie b^\perp}$. We may also replace
  $Z$ by $Z^B$ and $U^\CC$ by the analytic Zariski closure of $A$ in
  $U^\CC$ without changing our assumptions. With this in mind we have
  $\ker d\mua(x)=(\liea \cdot x)^\perp$ for all $x\in Z=Z^B$. This
  implies that $d\mu_{\lie b^\perp}(x)\colon T_x(Z^{(B)})\to \lie
  b^\perp$ is a bijection for all $x\in Z^{(B)}$. Since $\mu_{\lie
    b}\colon Z^B\to \lie b$ is locally constant on $Z^B$ this shows
  (3) and also (2).
\end{proof}

Let $A^c$ be the analytic Zariski closure of $A$ in $U^\CC$ and
$\liea^c$ its Lie algebra.  The group $A^c$ is a complex reductive
compatible subgroup of $U^\CC$ with maximal compact subgroup
$T=A^c\cap U$. We have $A^c=T\exp(\i\lie t)$ where $\lie t$ denotes
the Lie algebra of $T$.  In the following $\overline{S}$ denotes the
topological closure of a subset $S$ of $Z$ in $Z$. The same notation
is used for a subset of $\liea$ or more generally for subsets of a given 
topological space.

\begin{lemma} \label{lemma:strataconvexity} Let $S$ be an
  $\liea$-stratum of type $\liea/\lie b$ in $Z$. Then
  \begin{enumerate}
  \item $S$ is an $A^c$-stable locally closed complex submanifold of
    $Z$
  \item $\mua(\overline{S})=\overline{\mua(S)}$ is a convex polytope.
  \item every $y$ in $\overline S$ which is mapped by $\mua$ onto a
    vertex of $\mua(\overline S)$ is an $A$-fixed point.
  \end{enumerate}
\end{lemma}
\begin{proof}
  Since $A^c$ is connected and $A_{g \cdot y}=(A^c)_{g \cdot y}\cap
  A=(A^c)_y\cap A=A_y$ holds for all $g\in A^c$ and $y\in Z$ we have
  $(1)$.

Note that any $A$-stratum $S$ is an $A^c$-stable K\"ahlerian submanifold of
 $Z$. Let $\mu_{\lie t}\colon Z \to \lie t^*$ denote the momentum map on
 $Z$ given by restricting $\mu\colon Z\to \lieu^*$ to $\lie t$. In
 \cite{HeinznerHuckleberry} it is shown that  $\overline{\mu_{\lie t}(S)}$ is a
 convex polytope in $\lie t^*$. Equivalently $\overline{\mu_{\i \lie
      t}(S)}  = \mu_{\i \lie t}(\overline{S})$ is a convex polytope, where
 $\mu_{\i\lie t}\colon Z\to \i\lie t$ is the $A^c$-gradient map given by
$\mu$.  Since $\mua$ is the composition of $\mu_{\i \lie t}$ and  the orthogonal projection
  of $\i\lie t$ onto $\liea$ this shows that $\mua(\overline S)$ is a
  convex polytope in  $\lie a$.  Finally it is shown in
  \cite{HeinznerSchwarzStoetzel} that every $y\in S$ whose image is a
  vertex of $\mua(\overline S)$ has to be an $A$-fixed point.
\end{proof}

For the following Lemma we recall that we assume the $G$-action and
therefore also the $A$-action on $Z$ to be effective.

\begin{lemma}
  \begin{enumerate}
  \item There are only finitely many $A$-strata.
  \item The $A$-stable subset of $Z$ where $A$ acts freely is the
    unique open $A$-stratum is given by $S_0=\{z\in Z:\lie
    a_z=\{0\}\}$ and is open and dense in $Z$.
  \item $Z$ is the disjoint union of $A$-strata.
  \item The boundary $\overline S\setminus S$ of an $A$-stratum is a
    finite union of $A$-strata $\tilde S$ such that $\dim \tilde
    S<\dim S$ holds.
  \end{enumerate}
\end{lemma}

\begin{proof}
This follows from compactness of $Z$ and the Slice Theorem. 
\end{proof}

\begin{lemma}\label{lemma:isotropygroupofstrata}
  Let $S\not=S_0$ be an $\liea$-stratum of type $\liea/\lie b$ and
  $y\in S$.  Then there are $\liea$-strata $S_j$ of type
  $\liea/\liea_j$, $j=1,\dotsc,r$ such that $y\in \overline{S_j}$,
  $\dim\liea_j=1$ and $\lie b=\liea_1+ \dotsb + \liea_r$ hold.
\end{lemma}
\begin{proof}
  We fix a point $y\in S$ and apply the Slice Theorem to the
  $A$-action on $Z$ at $y$. This means that we find an $A$-stable open
  neighborhood $\Omega$ of $y$, a $B:=A_y$ subrepresentation $W$ of
  the isotropy representation $T_y(Z)$, an open $B$-stable
  neighborhood $\Omega_W$ of $0\in W$ and an $A$-equivariant
  diffeomorphism $\Psi\colon\twist AB\Omega_W\to \Omega$ such that
  $\Psi([e,0])=y$. Since the $A$-action on $Z$ is assumed to be
  effective and since it is real analytic the $B$-action on $W$ is
  effective. We view $\twist AB\Omega_W$ as an open subset of $\twist
  ABW$ and note that a $\lie b$-stratum $S(W)$ in $W$ of type $\lie
  b/\lie c$ determines uniquely the $\liea$-stratum $\twist
  ABS(W)\subset \twist ABW$ of type $\liea/\lie c$. This implies that
  we may restrict our attention to the $\lie b$ representation $W$.

  The image of $B$ in $\mathrm{GL}(W)$ is real diagonalizable since
  $B$ acts on $T_x(Z)$ by selfadjoint operators
  (\cite{HeinznerSchwarz}). Let $W=W_{\chi_0}\oplus \dotsb \oplus
  W_{\chi_r}$ be the isotypical decomposition of $W$ where for any
  linear function $\chi\colon \lie b\to\RR$ we set $W_\chi=\{w\in W:
  \xi\cdot w=\chi(\xi)w \text{ for all } \xi\in\lie b\}$ and $\chi_0$
  denotes the zero map.  We have $W_{\chi_0}=W^{\lie b}$. The open
  $\lie b$-stratum in $W$ is of type $\lie b$ and contains $W^{\lie
    b}\times (W_{\chi_1}\setminus\{0\})\times \dotsb\times
  (W_{\chi_r}\setminus\{0\})$. Then
  \begin{itemize}
  \item[a)] Any $\lie b$-stratum has $0$ in its closure and 
  \item[b)] if $0$ does not lie in the open $\lie b$-stratum, than
    there are $\lie b$-strata $S_j(W)$ of type $\lie b/\lie c_j$,
    $j=1,\dotsc, l$ such that $\dim \lie c_j=1$ and $\lie b=\lie
    c_1\oplus\dotsb\oplus \lie c_l$.
  \end{itemize}
  This follows from the fact that the $B$-representation $W$ is diagonalizable. 
  Since the $B$-action on $W$ is effective the open $\lie b$-stratum of $W$ is of type 
  $\lie b$.
\end{proof}

%%%%%%%%%%%%%%%%%%%%%%%%%%%%%%%%%%%%%%%%%%%%%
%%%%%%%%%%%%%%%%%%%%%%%%%%%%%%%%%%%%%%%%%%%%%
\section{Decomposition of the gradient map image}
\label{section:decomposition}
As in the previous section let $\liea$ be a linear subspace of $\liep$
which is a subalgebra of $\lieg$ and $A=\exp\liea$ the corresponding
commutative compatible subgroup of $G$.  Let $S$ be an $\liea$-stratum
of type $\liea/\lie b$. We set $\sigma:=\mu_\liea(S)$ and let
$\liea(S):=\liea(\sigma)$ be the unique affine subspace of $\liea$
which contains $\sigma$ as an open subset
(Lemma~\ref{lemma:strataconvexity}). We have 
$\liea(S)=\liea(\sigma)=q+\lie
b^\perp$ for any $q\in\mu_\liea(\overline S)$ where $\lie a=\lie
b\oplus \lie b^\perp$.  Since $\lie b$ only depends on $\liea(\sigma)$
we will also use the notation $\liea_\sigma=\lie b$. Formulated in more 
geometric terms $\liea_\sigma$ is the linear subspaces of $\liea$ which is
perpendicular to the affine linear space $\liea(\sigma)$ and  
coincides with the isotropy Lie algebra of any point $z\in S$.

Let $\Sigma:=\{\liea(\sigma): S\text{ is an $A$-stratum and
}\sigma=\mua(S)\}$ denote the set of all affine subspaces of $\liea$
obtained in this way.  For the open $A$-stratum $S_0$ we
have $\liea=\liea(\sigma_0)$ and $\sigma_0$ is the interior of
$P:=\mua(Z)$.  Let $\Sigma_1:=\{\sigma\in\Sigma:
\codim_{\liea}\liea(\sigma)=1\}$ and $P_0:=P\setminus
\bigcup_{\sigma\in\Sigma_1}P\cap \liea(\sigma)$.

\begin{lemma}\label{lemma:open}
The set $P_0$ is open in $\liea$. 
\end{lemma}
\begin{proof}
It is sufficient to show that every face $F$ of $P=\mua(Z)$
of codimension one is contained in $\liea(\sigma)$ for some 
$\sigma\in\Sigma_1$. 

The image $\mua(S)=\sigma$ of any $\liea$-stratum is open in
$\liea(\sigma)$. If we apply this to the open stratum $S_0$ we see
that for any face $F\not= P$ of $P$ this implies
that $S_0\cap\mua\inv(F)=\emptyset$. Since we have only finitely many
$\liea$-strata this shows for a face $F$ with $\codim_\liea F=1$ that
there is an $\liea$-stratum $S_F$ with $\sigma_F\in\Sigma_1$ and
that $\sigma_F$ is open in $F$.  We have $F\subset \liea(\sigma_F)$ and
therefore $P\setminus\bigcup_{\sigma_F}\liea(\sigma_F)$ is open in
$\liea$ where the union is over all faces of $P$ which are of
codimension one. This implies that $P_0$ is open in $\liea$.
\end{proof}

As in the previous section let $S_0$ denote the unique
open $\liea$-stratum in $Z$. 

\begin{lemma} \label{lemma:normalface}
We have $\mua\inv(P_0)\subset S_0$ or equivalently 
$\liea_y=\{0\}$ for all  $y\in\mua\inv(P_0)$. 
\end{lemma}
\begin{proof}
Assume that there is a $y\in \mua\inv(P_0)$ such that 
$\liea_y\not=0$. Let $\tilde S$ be the $\liea$-stratum 
which contains $y$. Since $\tilde S$ is not the open $\liea$-stratum
there is an $\liea$-stratum $S$ of type $\liea/\liea_1$ where $\dim\liea_1=\dim \liea-1$ 
such that $\tilde S\subset \overline S$. This shows that 
$\mua(y)\in\overline{\sigma}\subset\liea(\sigma)$ for $\sigma = \mu_\liea (S)$. Since 
$\sigma \in\Sigma_1$ this contradicts the definition
of $P_0$.
\end{proof}

Let $C(P_0)$ denote the set of connected components of $P_0$. 
For $\gamma\in C(P_0)$ let $P(\gamma)$ be the closure of the 
connected component $\gamma$.  The set $P(\gamma)$ is a convex polytope 
with non-empty interior $\int_\liea(P(\gamma))$ in $\liea$. Let 
$\mathcal{F}(P_0):=\{F: F \text{ is a face of } P(\gamma)
\text{ where } \gamma\in C(P_0) \}$ be the set of faces which are
determined by $P_0$. 
We have \( P_0=\bigcup_{\gamma\in C(P_0)}\int_{\liea}(P(\gamma)) \) and
\(P=\bigcup_{\gamma\in C(P_0)}P(\gamma)\). More importantly every face
$F\in \mathcal{F}(P_0)$ of codimension one is given by $P(\gamma)\cap
\liea(\sigma)$ for some $\sigma \in \Sigma_1$ and $\gamma\in C(P_0)$.

For a convex polytope $F$ in $\liea$ we introduce the following
notation. The affine span of $F$ is denoted by $\liea(F)$ and
$\int(F)=\int_{\liea(F)}(F)$ denotes the interior of $F$ as a subspace
of $\liea(F)$. The linear subspace of $\liea$ which is perpendicular to
$\liea(F)$ is denoted by $\liea_F$.  The dimension of $F$ is denoted
by $\dim F$ as is by definition the dimension of $\liea(F)$. Similarly
$\codim F$ means the codimension of $\liea(F)$ as a subspace of
$\liea$ and coincides with the dimension of $\liea_F$.
 
For $\gamma\in C(P_0)$ let $\Sigma_1(\gamma)$ denote the set
of codimension one faces of $P(\gamma)$.

\begin{proposition} \label{proposition:isotropy} Let $\gamma\in C(P_0)$
and let $F$ be a face of
  $P(\gamma)$ of codimension $k$. Then there are
  $\sigma_1,\dotsc,\sigma_k\in\Sigma_1(\gamma)$ such that
  \begin{enumerate}
  \item $F=P(\gamma)\cap \liea(\sigma_1)\cap\dotsb\cap
    \liea(\sigma_k)$ and
  \item \( \liea_y\subset
    \liea_F=\liea_{\sigma_1}+\dotsb+\liea_{\sigma_k} \) for all $q\in
    \int(F)$ and $y\in \mua\inv(q)$.
  \end{enumerate}
\end{proposition}

\begin{proof}
  Property (1) follows from the definition of $P(\gamma)$. 

  Let $q\in \int(F)$ and $y\in \mua\inv(q)$. Since
  $\sigma_1,\dotsc,\sigma_k\in\Sigma_1(\gamma)$ we have
  $\liea(F)=\liea(\sigma_1)\cap \dotsb\cap \liea(\sigma_k)$.  We have
  to show that $\liea_y\subset \liea_{\sigma_1}+ \dotsb+
  \liea_{\sigma_k}$.  Since $q\in \int(F)$ we have
  $\sum_{\sigma\in\Sigma_1 , q\in\liea(\sigma)}\liea_{\sigma}
  =\liea_{\sigma_1}+\dotsb+\liea_{\sigma_k}$.  The Slice Theorem
  implies that $ \liea_y=\sum_{\tilde
    \sigma\in\tilde\Sigma}\liea_{\tilde\sigma}$ for some subset
  $\tilde \Sigma\subset\Sigma_1$
  (Lemma~\ref{lemma:isotropygroupofstrata}).  This gives
  $\liea_y\subset \liea_{\sigma_1}+\dotsb+\liea_{\sigma_k}$.
\end{proof}

By the construction of the polytopes $P(\gamma)$ the set $\mathcal F(P_0)$ is closed under intersection and we have the following

\begin{remark} \label{non convex neighborhood(r)}
Let $D$ be the finite union of elements in $\mathcal F(P_0)$. Then the set of all points $\xi \in D$ such that $D$ is non convex in any neighborhood of $\xi$ is again a finite union of elements in $\mathcal F(P_0)$. 
\end{remark}

\section{Semistable points and convexity}
\label{semistable}
In this section we show that convexity of $A(Y)_+$ is closely related to
the behavior of semistable points after shifting. 

Let $\beta$ be a point in $\liep$. The $U$-orbit $U \cdot
\beta\subset\im\lieu$ can be identified with the coadjoint orbit $U
\cdot \im\beta \subset \lie u$ and is a complex flag manifold
$O:=U^\CC/Q$, where $ Q:= \{ g \in U^\CC \sodass  \lim_{t \to - \infty}
\exp(t\beta) \cdot  g \cdot  \exp(-t\beta) \ \text{exists in } U^\CC \}.  
$ In
particular, this induces a K\"ahler structure and  a holomorphic 
$U^\CC$-action on $U\cdot
\beta$. We denote this action  by $(g,x) \mapsto g\acts x$. The 
$G$-gradient map on $U^\CC \acts \beta$ is then just given by the  
projection of $O=U\cdot\beta\subset \i\lieu$ onto $\lie p$. 

\begin{proposition} {\rm (\cite{HSt05})}
  \label{Gbeta = Kbeta(p)}
  For $\beta \in \liep$ we have $G \acts \beta = K \cdot \beta$ in
  $O$.
\end{proposition}
The  $G$-gradient map $\mu_{\liep,\beta} : Z \times U^\CC \acts \beta \to \lie p,
(z,\xi) \mapsto \mup(z) - \pi_\liep(\xi) $ is called the shifting of $\mup$ with
respect to $\beta$ where $\pi_\liep\colon\im\lieu\to \liep$ denotes the orthogonal projection.  
In particular, $\beta$ is contained in the image
of $\mup\colon Z\to \liep$ if and only if 0 is contained in the image of
$\mu_{\liep,\beta}\colon Z\times U^\cc\acts \beta\to\liep$. 
The set of semistable points in $Y \times G\acts 
\beta \subset Y \times U^\CC \acts\beta$ with respect to the value $\alpha$ 
is by definition the set
\[
\mathcal{S}_G(\m_{\liep,\alpha})(Y\times G\acts\beta) := \{ (y, \xi)
\in Y \times G \acts \beta \sodass \overline{G \acts (y,\xi)} \cap
(\mu_{\liep,\beta})^{-1} (\alpha) \neq \emptyset \}
\]
for any $\alpha\in\liea$.  For $\alpha=0$ we
set $\mathcal{S}_G(\m_{\liep,0})(Y\times G\acts \beta)=
\mathcal{S}_G(\m_{\liep})(Y\times G\acts \beta)$. With this notation we have the following.

\begin{theorem} \label{convexity theorem(t)} Let $Y$ be a closed
  $G$-stable subset of $Z$ such that the intersection
  \[
  \mathcal{S}_G(\m_{\liep,\alpha_1})(Y \times G \acts \beta) \, \cap \,
  \mathcal{S}_G(\m_{\liep,\alpha_2})(Y \times G \acts \beta)
  \]
is nonempty for any $\alpha_j\in A_+(Y)$ and $\beta \in \liea$ with 
$\mathcal{S}_G(\m_{\liep,\alpha_j})(Y \times G \acts \beta)\not=\emptyset$.
Then $A_+(Y)$ is a convex polytope.
\end{theorem}

For the proof of the theorem we need some preparation. 
 
\begin{lemma}\label{lemma:kostant}
  Let $\liea_+$ be a closed Weyl-chamber in $\liea$ and
  $q,p\in\liea_+$. Then
  \[
  \lVert k\cdot q-p\rVert^2\ge\lVert q-p\rVert^2
  \]
  holds for all $k\in K$.
\end{lemma}
\begin{proof} Since the inner product on $\liep$ is $K$-invariant we have
\[
\lVert k\cdot q-p\rVert^2 - \lVert q-p\rVert^2 = -2 \, \cdot  <k \cdot q-q,p>.
 \]
We have $ <k \cdot q,p> 
= <\pi_\liea(k \cdot q),p>$ where $\pi_\ma$ is the orthogonal projection 
of $\lie p$ onto $\liea$. But $\pi_\ma(k \cdot q)$ 
is contained in the convex hull of the orbit of the Weyl group $W = N_K(\ma)/Z_K(\ma)$
through $q$ (\cite{Kostant}) and therefore it suffices 
to prove the inequality $<w \cdot q-q,p> \, \leq 0$ 
for all $w \in W$. Let $\alpha \in \Sigma^+$ and 
$\sigma_\alpha$ be the corresponding simple reflection. 
Then we have $<q,\alpha> \, \geq 0$ and therefore 
$<q,-\alpha> \, = \,  <q, \sigma_{\alpha} \cdot \alpha>\,  
= \, < \sigma_{\alpha} \cdot q,\alpha > \, \leq 0$ 
which implies $\sigma_\alpha \cdot q-q = - \lambda \alpha$ 
for some positive number $\lambda$. Since every Weyl group 
element can be written as a product of these simple 
reflections $w \cdot  q-q$ is a  negative linear 
combination of positive roots which shows $<w \cdot q-q,p> 
\, \leq 0$ for all $w \in W$.
\end{proof}

\begin{remark}
\begin{enumerate}
\item If one reads through our paper in the case that $G=A$, then one 
obtains a proof of the result we needed from Kostant's paper 
\cite{Kostant}. Thus our results are independent of \cite{Kostant}.
\item Kostant's fundamental paper  was the  first paper 
containing convexity results in the spirit presented here. 
\end{enumerate}
\end{remark}

\begin{proposition} \label{proposition:localminimumina} Let
  $p_0\in\liea$ and assume that $q_0 \in A(Y):=\mup(Y) \cap \liea$ is
  a minimum of the function $ \psi_{p_0} \colon A(Y)\to \RR$, $q
  \mapsto \lVert q-p_0 \rVert $.  Then
  \[
  (\mup\vert_Y)\inv(q_0)\subset Y^\xi := \{y \in Y \sodass \xi_Z(y)=0\}
  \]
  for $\xi:=q_0-p_0$.
\end{proposition}
\begin{proof}
We claim that $q_0$ is also a minimum of the function 
$\tilde\psi_{p_0}\colon \mup(Y) \to \RR $,  $q \mapsto \rVert q-p_0 \rVert$. 
Let $\liea_+$ be a Weyl chamber such that $p_0\in\liea_+$. 
Since $\mup(Y)$ is $K$-stable we have 
$\mup(Y)=K\cdot A(Y)=K\cdot (A(Y)\cap \liea_+)$. Let $\tilde q\in\mup(Y)$ 
and $k\in K$ such that $\tilde q=k\cdot q$ where $q\in \liea_+$.
This implies $\lVert\tilde q-p_0\rVert\ge\lVert q-p_0\lVert\ge\lVert q_0-p_0\rVert$
(Lemma~\ref{lemma:kostant}). 

Let $y\in (\mup\vert_Y)\inv(q_0)$ and $\xi=q_0-p_0$. Then $y$ is a
critical point of the function $\eta\colon G\cdot y\to \RR$, $g\cdot
y\mapsto \frac{1}{2}\lVert \mup(g\cdot y)-q_0\lVert^2$.  Now
$0=d\eta(y)=\langle d\mup(y),\mup(y)-p_0\rangle=d\mu^\xi(y)$ implies
$\xi_Z(y)=0$.
\end{proof}

\begin{proposition}\label{proposition:faces}
  We have $A(Y)\cap F=F$ for all faces $F\in\mathcal{F}(P_0)$ such
  that $\int(F) \cap A(Y)\not=\emptyset$.
\end{proposition}

We will prove the proposition recursively by arguing by dimension of
the faces of $\mathcal{F}(P_0)$ and starting with those faces which
are of maximal dimension. In order to carry this out we note the
following

\begin{lemma}\label{lemma:closureoffaces}
  Let $F^*\in \mathcal{F}(P_0)$. Then $F^*\cap A(Y)=F^*$ implies that
  $F\cap A(Y)=F$ for all $F\in \mathcal{F}(P_0)$ which are contained
  in $F^*$.\qed
\end{lemma}
\noindent
\textit{Proof of Proposition~\ref{proposition:faces}.}  
Let $F$ be an arbitrary face such that $\int_{\liea(F)}(F)\cap
A(Y)\not=\emptyset$. By the above indicated induction 
we may assume that our Proposition holds for all
faces $F^*\in\mathcal{F}(P_0)$ with $\dim F^*>\dim F$.
Lemma~\ref{lemma:closureoffaces} implies that we additionally may
assume that $\int(F^*)\cap A(Y)=\emptyset$ for all faces $F^*$ which
properly contain our given face $F$. The advantage of this assumption
is that for any $q_1\in \int(F)\cap A(Y)$ we can find a $r>0$ such
that $A(Y)\cap \Delta_r(q_1)=\int(F)\cap A(Y)\cap\Delta_r(q_1)$ holds.
For any $p_1\in\Delta_{\frac{r}{2}}(q_1)$ such that
$\xi_1:=p_1-q_1\in \liea_F$ is perpendicular to our face $F$ we obtain
$\lVert p_1-q_1\rVert\le\lVert p_1- q\rVert$ for all
$q\in\Delta_r(q_1)\cap A(Y)$.
Proposition~\ref{proposition:localminimumina} shows that
$(\xi_1)_Z(y)=0$ for all such $\xi_1$ and $y\in\mua\inv(q_1)$. Since
$q_1\in \int(F)\cap A(Y)$ was arbitrary this shows
$\liea_F\subset\liea_y$ for all $y\in \mua\inv(\int(F)\cap A(Y))$. Now
Proposition~\ref{proposition:isotropy} implies $\liea_y=\liea_F$ for
all $y\in \mua\inv(\int(F)\cap A(Y))$.

We will now argue that $\liea_y=\liea_F$ together with the assumption
that $\int(F)\cap A(Y)\not=\emptyset$ and $F\cap A(Y)\not=F$ leads to
a contradiction. 

Assume that there is a $q_1\in
\bd_{\int(F)}(\int(F)\cap A(Y)):=(\int(F)\cap
A(Y))\setminus\int_{\liea(F)}(F\cap A(Y))$ and let $r>0$ such that
$\Delta_r(q_1)\cap F\subset \int(F)$ and $\Delta_r(q_1)\cap A(Y)=F\cap
\Delta_r(q_1)\cap A(Y)$ hold.  Here we use the assumption that $
\int(F^*)\cap A(Y)=\emptyset$ for all faces $F^*$ which properly
contain $F$. Then there is a $p_0\in \Delta_{\frac{r}{2}}(q_1)\cap F$
with $p_0\not\in A(Y)$ and therefore a $q_0\in
A(Y)\cap\Delta_r(q_1)$ with satisfies $\lVert p_0-q_0\rVert\le\rVert
p_0-q\rVert$ for all $q\in A(Y)\cap\int(F)\cap\Delta_r(q_1)$.  Since
$A(Y)\cap\int(F)\cap\Delta_r(q_1)= A(Y)\cap\Delta_r(q_1)$ we may 
apply Proposition~\ref{proposition:localminimumina}.
This gives $\xi_Z(y)=0$ where $\xi=p_0-q_0$ and $y\in \mua\inv(q_0)$.
Since $\xi\not=0$ and $\xi\not\in \liea_Y$ this contradicts
$\liea_F=\liea_y$.  \qed

\begin{corollary}\label{corollary:maincommutative}
  The set $A(Y)$ is a union of faces $F\in \mathcal{F}(P_0)$
and is therefore a finite union of convex polytopes each of it the 
convex hull of images of fixed points of $T$ in $Y$.\qed
\end{corollary}

For the proof of Theorem~\ref{convexity theorem(t)} we also need the
the fact that a subset $D$ of an
  Euclidian vector space which is a finite union of convex polytopes
  and is not convex has the property that for any sufficiently small $r > 0$ there
  exists a point $\beta \in \liea$ such that the closed ball of radius
  $r$ and center $\beta$ meets $D$ in precisely two points $\alpha_1$
  and $\alpha_2$. This geometric input has also been used in Kirwan`s proof
of her convexity result (\cite{Kirwan2}).
\bigskip

\noindent {\it Proof of Theorem}~\ref{convexity theorem(t)}. The set
$A_+(Y)$ is a finite union of convex polytopes
(Corollary~\ref{corollary:maincommutative}).  Assume that $A_+(Y)$ is
not convex. Then there exist $r > 0$ and $\beta \in
\liea$ such that the closed ball of radius $r$ and center $\beta$
meets $A_+(Y)$ in precisely two points $\alpha_1$ and $\alpha_2$ which 
are on the boundary of this ball.  Now if for $\alpha\in\liea$ the value
$\lVert\alpha\rVert$ is critical for the function
$\lVert\mu_{\liep,\beta}\rVert^2\colon Z\times G\acts\beta\to\RR$ then
there is an associated pre-stratum $S_\alpha=\{w\in Z\times G\acts
\beta: \lVert\alpha\rVert= \min\{\lVert\mu_{\liep,\beta}(g\cdot
w)\rVert: g\in G\}\}$ for the $G$-action on $Z\times G\acts \beta$ in
the sense of \cite{HeinznerSchwarzStoetzel}.  The values
$\lVert\alpha_j\rVert$ are critical points of
$\lVert\mu_{\liep,\beta}\rVert^2\colon Z\times G\acts\beta\to\RR$ and
define two non empty $G$-pre-strata for the $G$-action on $Z\times
G\acts \beta$.  We have $S_{\alpha_j}\cap Y=
\mathcal{S}_G(\m_{\liep,\alpha_j})(Y\times G\acts\beta)$.  Since
$\alpha_1\not=\alpha_2$ and $\alpha_j\in\liea_+$ these pre-strata are
disjoint. Consequently we have
\[
\mathcal{S}_G(\m_{\liep,\alpha_1})(Y \times G \acts \beta) \ \cap \
\mathcal{S}_G(\m_{\liep,\alpha_2})(Y \times G\acts \beta) =
\emptyset\,.
\]
This contradicts the assumption of Theorem~\ref{convexity theorem(t)}.
\qed

\begin{remark}
   In the case where $G=U^\CC$ and $Y$ is an irreducible $U^\CC$-stable
  complex subspace of $Z$ the set $Y\times U^\CC \acts \beta$ is a
  K\"ahlerian space and the pre-strata in the above proof are locally
  closed complex subspaces of $Y\times U^\CC\acts \beta$.  This
  implies that there is a unique open $U^\CC$-stratum which is dense
  in $Z\times U^\CC \acts \beta$. The above proof then gives Kirwan's
  convexity theorem for actions of complex reductive groups on compact
  $U^\CC$-stable irreducible complex subspaces of K\"ahler
  manifolds. This is a rather special case of the more general
  convexity result in \cite{HeinznerHuckleberry}.
\end{remark}

\section{The projective case}\label{projective}

We fix now a finite dimensional unitary
representation space $V$ of the compact group $U$ and consider
$Z=\mathbb{P}(V)$.  The action of $U$ on $V$ extends to a holomorphic
linear action of $U^\CC$ and induces an algebraic $U^\CC$-action on
the associated complex projective space $\mathbb{P}(V)$.  There are
$G$-gradient maps $\mu_{\liep,V}\colon V\to \liep$,
$\mu_{\liep,V}^\xi(v)=\langle \xi v,v\rangle$ on $V$ with respect to
the K\"ahler structure induced by the unitary one on $V$ and a
$G$-gradient map $\mu_{\liep,\mathbb{P}(V)}\colon \mathbb{P}(V)\to \liep$,
$\mu_{\liep,\mathbb{P}(V)}^\xi([v])=\frac{\langle \xi v,v\rangle}{\lVert
  v \rVert^2}$ on $\PP(V)$ with respect to the induced Fubini-Study
K\"ahlerian structure on $\PP(V)$. Here we denote the fixed positive
Hermitian structure on $V$ by $\langle\ ,\ \rangle$ and $[v]\in\PP(V)$ 
denotes the line through $v\in V\setminus\{0\}$. Note that the Fubini-Study
form on $\PP(V)$ is given by symplectic reduction and is up to a positive
constant the unique K\"ahler form on $\PP(V)$ which is invariant with respect to the natural 
action of the special unitary group $\mathrm{SU}(V)$. 

In order to simplicity the notation we set $\mup:=\mu_{\liep,\PP(V)}$.
We view $\PP(V)$ as a real algebraic variety and fix a $G$-stable closed real
semialgebraic subset $Y$ of $\PP(V)$.  We say that $Y$ is irreducible
if the real Zariski closure of $Y$ in $\PP(V)$ is a real irreducible
subvariety. Our main result is 
 
\begin{theorem} \label{convexity projective(t)} The set   $A_+(Y) := \mup(Y) \cap \ma_+$ is a convex polytope.
\end{theorem}

We have the following consequences which are shown below.

\begin{corollary} \label{flag convexity}
Let $Z = U^\cc/Q$ be a complex flag manifold with $G$-gradient map 
$\mu_\liep :Z \to \liep$. Then the sets
$A_+(Z)$ and $A_+(\overline{G \cdot x})$
are convex polytopes.
\end{corollary}

Using this fact we also have

\begin{corollary} \label{open orbits}
Let $Z = U ^\cc/Q$ be a complex flag manifold and assume that  
$G$ is a real form of $U^\cc$ which is given as the set of fixed points
of an anti-holomorphic involution commuting with the given Cartan involution 
on $U^\cc$. Let $\xi\in A_+(Z)$ be the unique closest point to the origin.
The set of semistable points  $\mathcal S_G(\mathcal M_{\lie p,\xi})(Z) 
:= \{z \in Z \sodass \overline{G \cdot z} \cap \mup\inv(\xi) \neq \emptyset\}
$ coincides with the union of all open $G$-orbits in $Z$. Moreover, the closed 
$K^\cc$-orbits have the same image under the $G$-gradient  map  $\mip$.
\end{corollary}

To prove Theorem \ref{convexity projective(t)} we use the same strategy as in the proof of Theorem~\ref{convexity theorem(t)}. For this we need the notion of quasi-rational points in $\lie a$. 
Let $\lie h$ be a maximal abelian subalgebra of the centralizer $\lie z_{\lie k}(\lie a)$. Then $\lie s_\mun := \lie h \oplus \im \lie a$ is maximal torus in $\lie u$. We call a point $\alpha \in \lie a \simeq \im \lie a$ quasi-integral if $\alpha$ is the projection of an integral element $\alpha^\prime$ in the compact torus $\ms_\mun = \mh \oplus \im \ma$ onto $\im \lie a$. We denote this projection by $\pi_{\im \lie a}$.  A point $\beta \in \im \lie a$ is  called quasi-rational if it is a rational multiple of an quasi-integral element $\alpha \in  \ma$.
The following lemma allows us to
make a  reduction to quasi-rational points.

\begin{lemma} \label{rational(l)}
The set $A(Y)_+$ is a finite union of 
quasi-rational polytopes, i.\,e. the polytopes are convex hulls of finitely many quasi-rational points. In particular, the quasi-rational points are dense in $A(Y)_+$.
\end{lemma}

\begin{proof}
Corollary \ref{corollary:maincommutative} says that $A(Y)_+$ is 
the intersection of a positive Weyl chamber with a finite union of 
convex polytopes which are given by the
convex hull of images of sets of $A$-fixed points in $Y$. 
If $[v] \in Y^A$. Then $v$ is 
contained in a weight space
$
V_\chi := \{ v \in V \sodass \xi\cdot v = \chi(\xi) \cdot v 
\ \forall \  \xi \in \liea\}
$
of the $\liea$-representation $V$. Here $\chi \colon \liea \to \RR$ is
a linear function with $\chi = \im \varphi_\ast\vert_{\im \lie a}$ for some character $\varphi \colon S_U \to S^1$. Here $S_U$ denotes the maximal torus of $U$ with Lie algebra $\lie s_\lie u$ and $S^1$ is the maximal compact subgroup of $\CC^*$.  For every $\xi \in  \ma$ we therefore have 
$
\mia([v])(\xi) = \chi(\xi).
$
So $\mia([v])$ is an  quasi-integral element in $\ma$ in the sense of the appendix. Consequently $A(Y)_+$ 
is a finite union of quasi-rational convex polytopes.
\end{proof}

Any semialgebraic set $Y$ has a finite semialgebraic stratification,
i.\,e.  $Y$ can be decomposed into real analytic locally closed
semialgebraic submanifolds $A_i$ such that for $\overline{A_i} \cap
A_j \neq \emptyset$ we have $A_j \subset \overline{A_i}$ and $\dim A_j
< \dim A_i$.  We choose one decomposition, define $\hat Y$ to be the
union of the maximal dimensional strata and set $\dim Y := \dim \hat
Y$.  If $Y$ is a real algebraic set this gives the Krull dimension of
$Y$.  In general we have $\dim Y = \dim cl(Y)$ where $cl(Y)$ denotes
the real Zariski closure of $Y$. For detailed proofs see e.\,g.
\cite{BenedettiRisler} or \cite{Coste}. 

Let $Y$ be a
  closed $G$-stable irreducible semialgebraic subset of $\PP(V)$.
The technical part of the proof of Theorem~\ref{convexity projective(t)}
is the following

\begin{proposition} \label{SGMipalpha offen und dicht(p)}  Let
  $\alpha \in \ma_+$ and $\beta\in \liea$ such that 
  \[
  \mathcal S_{G} (\m_{\liep,\alpha})(Y\times G\acts\beta)
  \]
  is non empty. Further assume that $\beta$ is quasi-rational. Then 
 \[
  \mathcal S_{G} (\m_{\liep,\alpha})(Y\times G\acts\beta) \cap
  (\hat Y\times G\acts\beta)
  \]
is open and dense in $\hat Y\times G\acts \beta$.
\end{proposition}

Actually the assumption that $\beta$ is quasi-rational is not necessary but simplifies the proof.

\medskip
\noindent
{\it Proof of Theorem~\ref{convexity projective(t)}}.
Assume $A_+(Y)$ is non convex. As in the proof of Theorem~\ref{convexity theorem(t)} we get points $\alpha_1, \alpha_2 \in \ma_+$ and $\beta \in \ma$ such that  
\[
\mathcal{S}_G(\m_{\liep,\alpha_1})(Y \times G \acts \beta) \ \cap \
\mathcal{S}_G(\m_{\liep,\alpha_2})(Y \times G\acts \beta) =
\emptyset\,.
\]
Using Lemma \ref{rational(l)} and Remark \ref{non convex neighborhood(r)} one can choose the point $\beta$ to be quasi-rational. See \cite{Kirwan2} for the explicit construction.
But this contradicts Proposition~\ref{SGMipalpha offen und dicht(p)} and shows the assertion.

\qed
\bigskip

\noindent
{\it Proof of Proposition~\ref{SGMipalpha offen und dicht(p)}.} 
Since $\beta$ is quasi-rational, there exists an integral element $\delta^\prime \subset \lie s_\lie u$ such that $\beta = \frac \im n \cdot \pi_{\im \lie a}(\delta^\prime)$ 
for some $n \in \nn$. Therefore we have a $G$-equivariant 
diffeomorphism
$
\varphi \colon  Y \times G \acts \beta^\prime \to Y^\prime \times G \acts \delta^\prime
$
given by multiplication with $n$ in the second component. Here $\beta^\prime = \frac 1 n \cdot  \delta^\prime$ and  $Y^\prime = Y$ 
but seen as a subset of $\pp(V)$ endowed with the Fubini Study form multiplied 
with $n$. In particular, we have
\[
\mathcal M_{\lie p,\alpha}(Y \times G \acts \beta) = 
\mathcal M_{\lie p,\alpha}(Y \times G \acts \beta^\prime) = 
\varphi\inv (\mathcal M_{\lie p, n \cdot \alpha}(Y^\prime \times G \acts \delta)).
\]
Since $\varphi$ is $G$-equivariant, we get the analog equation for the sets 
of semistable points. Since $Y^\prime \times G \acts \delta^\prime$ is a closed $G$-stable irreducible 
semialgebraic set of some projective space it suffices to prove the proposition for $\beta = 0$.

So let $\alpha \in \ma_+$ such that $\mathcal M_{\lie p,\alpha}(Y)$ is non empty. Then there exists a point $y \in Y$ such that $\mip(y) = \alpha$. Since the quasi-rational points are dense in $A_+(Y)$ there exists a sequence $\seq{y_n}$ such that $\mip(y_n) =: \alpha_n$ are quasi-rational and $\lim_{n \to \infty} \alpha_n = \alpha$. In particular $\mathcal M_{\lie p, \alpha_n}(Y)$ is non empty.

Let assume that the open set $\mathcal S_G(\Mipa)(Y) \cap \hat Y$ is not dense in $\hat Y$. Then there exists a $G$-stable open subset $U$ in the complement and an $r >0$ such the $\Vert \mip(y^\prime) - \alpha\Vert^2 \geq r$ for all $y^\prime \in U$. Moreover, there exists an $N \in \nn$ such that for all $n \geq N$ we have 
$\Vert \mip(y^\prime) - \alpha_n\Vert^2 \geq \frac r 2$ for all $y^\prime$ in the $G$-stable subset $U$. Therefore, for $n \geq N$, the set 
$\mathcal S_G(\mathcal M_{\lie p, \alpha_n})(Y) \cap \hat Y$ is a non empty open subset of $\hat Y$ which is not dense.  Therefore it suffices to prove the assertion for quasi-rational $\alpha$.

If $\alpha \in \ma_+$ is a quasi-rational point, i.\,e. 
$\alpha = \frac p q \cdot \gamma \in \ma_+$ for some quasi-integral element 
$\gamma$ and some coprime integers $p$ and $q$, it follows 
>from the results given in the second part of the appendix that
$
Y_\alpha = \Phi_\alpha(Y \times K \cdot [v_\gamma^\prime])
$
is a closed $G$-stable irreducible semialgebraic subset of $\pp(V_\alpha)$. Here $\gamma^\prime$ again is a integral element in $\lie s_\lie u$ over $\gamma$. 
Note that the projection of $\Phi_\alpha\inv(\SGM(Y_\alpha) \cap \hat Y_\alpha) \subset Y \times [v_{\gamma^\prime}]$ is just $\mathcal S_{G} (\mathcal M_\lie p(\alpha))(Y) \cap \hat Y$ where $\hat Y_\alpha = \Phi_\alpha(\hat Y \times K \cdot [v_{\gamma^\prime}])$. Therefore we can also restrict to the case $\alpha= 0$.
 
 We 
have $\mp(Y)\not=\emptyset$ and 
$
\SGM (Y)= Y \backslash(Y \, \cap \, \pi(\zcone))
$ where $\mathcal N_G := \{ v \in V \sodass 0 \in \overline{G \cdot v} \}$ is the null cone in $V$. 
Since the null cone is a real algebraic subset of 
$V$ (Lemma~\ref{algebraicity}) and $cl(Y)$ is irreducible, 
the intersection $\pi(\zcone) \cap cl(Y)$  is either $cl(Y)$ or a proper 
algebraic subset of lower dimension in $cl(Y)$. 
Since $\dim \hat Y = \dim Y = \dim cl(Y)$
and $\Mip(Y) \neq \emptyset$ the set  $\hat Y \cap \pi(\zcone)$ is a proper 
semialgebraic subset of lower dimension in $\hat Y$. In particular, its 
complement $ \SGMip Y \cap \hat Y$ is open and dense in $\hat Y$ 
(see e.\,g. \cite{BenedettiRisler} for basic properties of semialgebraic sets).

\qed

\medskip
We now prove the two corollaries \ref{flag convexity} and \ref{open orbits} 
about complex flag manifolds. Note first that every complex flag manifold 
$Z = U^\cc/Q$ can be identified with an orbit $U \cdot \beta$ where $\beta$ 
is contained in the cone 
$C_Q := \{ \lambda \in \ms_\mun \sodass \lambda = \sum c_j  \cdot \alpha_j ,  
\alpha_j \in \Pi^\prime, c_j \in \rr^+\}$ where $\Pi^\prime$ is the subset 
of the simple roots which define $Q$.
Using the notation of the second part of the appendix we have the following

\begin{corollary}
Let $Z = U^\cc/Q$ be a complex flag manifold and let $Y$ be a closed $G$-stable 
subset such that $\varphi_\beta(Y)$ is an irreducible semialgebraic subset  of 
$\pp(\Gamma_{\beta})$ for every integral $\beta$ in the cone $C_Q$.
Then $A_+(Y)$ is a convex polytope for every $G$-gradient map 
$\mup \colon Z \to \lie p$ on $Z$.
\end{corollary}

\begin{proof}
Every momentum map on the complex flag manifold $Z$ with respect to the 
$U$-action is of the form
$
\mu \colon  Z \to \mun$, $ x \mapsto \alpha_x,
$
where $\alpha$ is an element in the cone $C_Q$. By Theorem 
\ref{convexity projective(t)}, the lemma holds if $\alpha$ is 
an integral element in $C_Q$. If $\alpha$ is a rational point, it can 
be written
in the form $\alpha = \frac 1n \cdot\beta$ for some integral element 
$\beta$ in $\ms_\mun$ and some $n \in \nn$.
In particular, $\mip(Y) = \frac 1n \cdot \tilde 
\mu_\lie p(Y),$
where $\tilde \mu$ is given by $\tilde \mu(x)= \beta_x$. 
This proves the rational case.
Since the rational points are dense in each cone 
$C_Q$, we can construct a sequence $\seq{\alpha_n}$ of rational 
points $\alpha_n \in \ms_\mun$ such that $Q_-(\alpha_n)$ coincides with
$Q$ for all $n \in \nn$ and
$
\lim_{n \to \infty}\alpha_n = \alpha.
$
This gives a sequence of momentum maps $\seq{\mu^{(n)}}$ such that 
the assertion holds for every element of this sequence. This gives the 
corollary for the limit point $\mu$.
\end{proof}
In particular, the above corollary can be applied to the complex 
flag manifold $Z$ itself and closures of $G$-orbits in $Z$ which 
gives Corollary \ref{flag convexity}.

\medskip
\noindent
{\it Proof of Corollary \ref{open orbits}}. \ 
Since $\xi$ is the unique closest point to the origin in 
$A_+(Z)$, a point $x_0 \in (\mip)\inv(\xi)$ is a global minimum of 
$\eta_\lie p$. Therefore, the $G$-orbit through $x_0$ is open and 
the $K^\cc$-orbit  through $x_0$ is closed and coincides with the $K$-orbit through $x_0$ (see \cite{Bremigan} and \cite{Vilonen}).
So we need to show that all open  $G$-orbits are contained in the 
set $\SGMip{\xi}$. By Proposition \ref{SGMipalpha offen und dicht(p)},
 we know that $\SGMip{\xi}$ is open and dense in $Z$ in the integral case. 
This can be extended to the non integral case as in the proof of 
Proposition \ref{SGMipalpha offen und dicht(p)}. Consequently, the 
set  $\SGMip{\xi}$ contains all open $G$-orbits and the closed $K^\cc$-orbits 
in these open $G$-orbits are contained in  $(\mip)\inv(K \cdot \xi)$.
\qed

\section{Appendix}

\subsection{Algebraicity of the null cone} In the proof of 
Proposition \ref{SGMipalpha offen und dicht(p)} we use the 
algebraicity of the null cone $\mathcal N_G= \{ v \in V \sodass 0 \in \overline{G \cdot v} \}$. 
Here we give a proof of this fact.
 
\begin{lemma} \label{algebraicity}
The null cone $\mathcal N_G $ is an algebraic subset of $V$.
\end{lemma}

\begin{proof}
Using the decomposition of $G$ into its semisimple part $G_s$ and its center $Z(G)$ 
(see \cite{HSt05}) we first prove the algebraicity of the null cone $\zconess$ with 
respect to $G_s$
by using results of \cite{Birkes} and \cite{RichardsonSlodowy}. 
All of them appear as special cases in \cite{HeinznerSchwarz}.

Let  $G_s^\cc$ be the complexification of $G_s$ and let $V^\cc$ 
the corresponding complexified representation space. Further let 
$F_1,...,F_d$ be the generators of the  algebra $\cc[V^\cc]^{G_{s}^\cc}$ 
of $G_s^\cc$-invariant polynomials on $V^\cc$. Then the map 
$
F = (F_1,..,F_d) \colon  V^\cc \to \cc^d,
$
parameterizes the Zariski closed $G_s^\cc$-orbits in $V^\cc$. 
Without loosing generality we can assume that $F_j(0)=0$ for all generators $F_j$.

By a result of \cite{Birkes}, for every $v \in V$ the orbit $G_s^\cc \cdot v$ is 
Zariski closed in $V^\cc$ if $G_s \cdot v$ is closed in $V$. Therefore, 
for every closed orbit $G_s \cdot v$, $v \neq 0$,
there exists a function $F_j$ from the list of generators of the algebra  
$\cc[V^\cc]^{G_s^\cc}$ such that
$
0 \neq F_j\vert_{G_s \cdot v}
$.
Since $G_s$ is Zariski dense in $G_s^\cc$,  we can assume that the 
polynomials $F_j$ are extensions of real polynomials $f_1,..,f_d$ 
which generate the algebra $\rr[V]^{G_s}$.
In particular, we have
$
0 \neq f_j\vert_{G_s \cdot v}
$
which shows that the null cone is given as the real algebraic subset
$
\{ v \in V \sodass f_j(v) = 0 , \ \ j = 1,...,d\}
$
of $V$.

To prove the general case let $V_j$ denote the weight spaces of $Z(G)$ in $V$. 
Since $G_s$ 
and $Z(G)$ commute, these subspaces are stable under $G_s$. In particular, 
we have
$
\rr[V]^{G_s} = \bigotimes \rr[V_j]^{G_s}.
$
Each factor $\rr[V_j]^{G_s}$ has finitely many generators which can be 
chosen to be homogeneous polynomials. 
Let $f\colon V \to \rr^k$ be the polynomial map which is 
given by all these generators. Then $f$ is invariant with respect 
to $G_s$ and equivariant with respect to $Z(G)$. Here the action 
of $Z(G)$ on $\rr^k$ is given by the action on each generator.
The corresponding null cone of $Z(G)$ in $\rr^k$ is a finite union of 
linear subspaces $H_j \subset \rr^k$ (see \cite{HeinznerSchwarz} Corollary 15.5). 
Therefore, the 
preimage of this null cone under $f$ is an algebraic subset of $V$ which we call 
$\mathcal N^\prime$. We show that the null cone $\zcone$ coincides with this 
algebraic set $\mathcal N^\prime$.

Let $G_U$ denote the maximal 
compact subgroup of $G_s^\cc$ and let $H$ be a $G_U$-invariant 
positive definite Hermitian form on $V^\cc$ such 
that the alternating part vanishes on $V$. Such a 
form exists by \cite{RichardsonSlodowy} and we get a momentum map
$\mu_{G_U} \colon  V^\cc \to \mg_U^\ast$
for the $G_U$ action on $V^\cc$. By construction,
$\tMip = \mathcal M \cap V$ where $\mathcal M := (\mu_{G_U})\inv(0)$. 

We have $\zcone \subset \mathcal N^\prime$. For the opposite inclusion 
let $v \in \mathcal N^\prime$. By definition, there exists a sequence 
$\seq{g_n} \subset Z(G)$ such that 
$
\lim_{n \to \infty} g_n \cdot f(v) = 0.
$
For every $g_n \cdot f(v) \in \rr^k$ let $\alpha_n$ be a point in a closed 
$G_s$-orbit in $f\inv(g_n \cdot f(v))$. By \cite{RichardsonSlodowy} every 
closed $G_s$-orbit intersects the set $\tMip$ and we may choose 
$\alpha_n \in \tMip$. Let $F$ be the complex extension of $f$ to 
$V^\cc$. The map $f$ is a proper map when restricted to $\mathcal M$ 
(see \cite{RichardsonSlodowy}). This is also true for the restriction 
of $F$ to $\tMip$ since
$\tMip = \mathcal M \cap V$. Consequently the sequence $\seq{\alpha_n}$ 
has a convergent subsequence with limit point $\alpha \in \tMip$. 
But $f(\alpha) =0$ which implies $\alpha \in \zconess$ and consequently 
$v \in \zcone$.
\end{proof}

\subsection{Shifting with respect to quasi-rational points}

Since there is in general no symplectic embedding of the orbit 
$U \cdot \beta \subset \im \mun$ into a projective space $\pp(W)$, 
equipped with the Fubini-Study metric coming from a $U$-invariant
 Hermitian form on $W$, the set $Y \times G \acts \beta$  
is in general not contained in the class of examples we are considering 
in section \ref{projective} of this paper. Therefore, we have to introduce 
a slight modification of the shifting procedure.

Given an integral element
$\alpha^\prime \in \ms_\mun \subset \ms$ we get an associated character 
$\chi_{\alpha^\prime} \colon Q \to \cc^\ast$ on the parabolic subgroup and  a
 $U^\cc$-homogeneous line bundle
\[
L^\alpha = U^\cc \times_{\chi_{\alpha^\prime}} \cc = (U^\cc \times \cc)/Q
\]
over $U^\cc/Q$. Let $\Gamma_{\alpha^\prime} := \Gamma(U^\cc/Q,L^{\alpha^\prime})$ denote 
the space of holomorphic sections of the line bundle $L^{\alpha^\prime}$ and 
let $\Gamma_{\alpha^\prime}^\ast$ denote its dual space.
By the theorem of Borel and Weil, the space $\Gamma_{\alpha^\prime}$ is an 
irreducible $U^\cc$-representation space with highest weight $\alpha^\prime$ and 
it follows that  
the projective space $\pp(\Gamma_{\alpha^\prime}^\ast)$ contains a unique complex 
$U$-orbit $U \cdot [v_{\alpha^\prime}]$ 
with $ U \cdot [v_{\alpha^\prime}]  \simeq U^\cc/Q
$. In particular this gives an embedding of $U \cdot \alpha^\prime$ into a projective space.
For detailed proofs see e.\,g. \cite{Akhiezer} or \cite{Huckleberry}.

Now, let $\beta \in  \ma$ be a quasi-rational point and let $\alpha$ be the 
unique minimal quasi-integral element in $\rr^+ \cdot \beta$. 
Then $\beta = \frac p q \cdot \alpha$
for coprime natural numbers $p$ and $q$. Let $\alpha^\prime$ be an integral element in $\lie s_\lie u$ with $\pi_{\im \lie a}(\alpha^\prime) = \alpha$.
 The momentum map 
$\mu_{\pp(\Gamma_{-\alpha^\prime}^\ast)}$ restricted to the orbit 
$U \cdot [v_{- \alpha^\prime}]$ is given by
$
\mu_{\pp(\Gamma_{- \alpha^\prime}^\ast)}([v])
= \mu_{\pp(\Gamma_{- \alpha^\prime}^\ast)}(u \cdot [v_\alpha^\prime])
= - u \cdot \alpha^\prime =: - \alpha^\prime_{[v]}.
$
Since $\beta$ and $\alpha$ are related by the fixed numbers $p$ and $q$, 
we can define a unique embedding
\[
\Phi_\beta \colon 
Y \times U \cdot \alpha^\prime \hookrightarrow \pp(V) \times \pp(\Gamma_{-\alpha^\prime}^\ast)
\hookrightarrow
\pp( V^{\otimes q}\otimes (\Gamma_{-\alpha^\prime}^\ast)^{\otimes p}) =: \pp(V_\beta)
\]
using the Segre embedding.
Define $Y_\beta := \Phi_\beta(Y \times G \acts  \alpha^\prime)$.
Since $G \acts \alpha^\prime$ is a real algebraic set we
may choose $\hat Y_\beta$ to be $\Phi_\beta( \hat Y \times G \acts \alpha^\prime)$. 
We get a $K$-equivariant map
\[
\mipb \colon  Y_\beta \to \lie p^\ast, \quad \Phi_\beta((y,\xi)) \mapsto
q \cdot \mip(y) - p \cdot \xi
\]
which we call the shifting of $\mip$ with respect to the quasi-rational point $\beta \in  \ma$. 
In particular, this is contained in the class of examples which we consider 
in section \ref{projective}.

\end{document}